\documentclass[]{article}
\usepackage{latexsym,amssymb,amsmath,amsfonts,epsf,rotating}
\usepackage[dvips]{epsfig}

\begin{document}

\vspace{5mm}

\centerline{\Large\bf Implications of a New Characterisation of }
\centerline{\Large\bf the Distribution of Twin Primes}

\vspace{5mm}

\centerline{P.F.~Kelly\footnote{patrick\_kelly@ndsu.nodak.edu}
 and Terry~Pilling\footnote{terry@mailaps.org}}

\vspace{5mm}

\centerline{Department of Physics}
\centerline{North Dakota State University}
\centerline{Fargo, ND, 58105-5566}
\centerline{U.S.A.}

\vspace{5mm}

\begin{abstract}
We bring to bear an empirical model of
the distribution of twin primes and produce two distinct results.
The first is that we can make a quantitative probabilistic prediction
of the occurrence of gaps in the sequence of twins within the primes.
The second is that the ``high jumper'' {\it i.e.,} the separation with
greatest likelihood (in terms of primes) is always expected to be zero.

\end{abstract}

\vspace{5mm}


\vspace{5mm}

\noindent
MCS:  11N05 (Primary) 11B05, 11A41 (Secondary)

\vspace{5mm}

\section{Introduction}

In a recent paper~\cite{random}, 
hereafter referred to as ``{\it TPI}'' we introduced a novel empirical 
characterisation of the sequence of twin primes.
Notation, nomenclature and conventions developed in 
{\it TPI} are carried over here.
There are three salient features of the model.
First is the viewpoint that the sequence of twins is most naturally 
and usefully considered as a subset of the primes rather than in terms 
of the natural numbers.
The second feature is that for the finite subset of primes less than
some number $N$, twins occur in the manner of fixed-probability 
random events.
{\it I.e.,} there is a constant probability that each successive prime
following a given twin is itself the first member of the next twin.
The third feature is that this ``constant'' probability is not 
universal.
It varies in a rather simple way with the {\it length} of the sequence of
primes up to $N$.

In this paper, we investigate a number of aspects of and make predictions
which follow directly and indirectly from this new model for the 
distribution of twins.
In the next section, we shall briefly describe the details of 
our empirical model.
The following section consists of two separate results.
First are probabilistic predictions for the occurrences of 
large gaps in the sequence of twins.
These predictions are verified by direct comparison
in the regime in which the model was constructed.
One important comment that must be made is that the gaps that we note
are separations in the sequence of primes.
Those usually quoted in the literature are {\it arithmetic differences},
or spacings in the natural numbers.
The second result we report is that the so-called ``high-jumper'' 
(the most probable prime separation) is zero in our model.
Our results do not admit a direct comparison with 
high-jumper analyses extant in the literature which are concerned
with most probable natural number separations between twins
less than a given $N$.

\section{The Model}

The model that we consider is empirical in that it is derived from
a direct analysis of the distribution of twin primes less than
$2 \times 10^{11}$.
The essential feature which provides the key to the success of the
model is that the distribution of twins is considered in the context
of the primes alone rather than within the natural numbers.
The model is based upon the observation that twins less than
some number $N$ seem to occur
as fixed-probability random events in the sequence of primes.
That is, there is a characteristic distribution of {\it prime separations}
which may expressed in the form
\begin{equation}
{\cal P}(s, \pi_1) = m \, e^{-ms} \, .
\label{prob}
\end{equation}
Here, $\pi_1$ is the number of primes less than or equal to $N$,
$s$ is the {\it prime separation} (the number of unpaired {\it singleton}
primes occuring between two twins),
and $m$ is a decay parameter which is constant for a given $N$,
but varies with $\pi_1$.
${\cal P}(s, \pi_1)$ is the continuous probability density that a given 
pair of twins in the sequence of primes up to $N$ has prime separation 
$s$.
It is properly normalized: $\int_0^\infty {\cal P}(s, \pi_1) \, ds \equiv 1$.

We chose a representative sample of prime sequences
and determined the decay constants for each.
We began our analysis with $(5 \ 7)$, 
discarding the anomalous twin $(3 \ 5)$.
The variation of the decay parameters -- the slopes on a plot of 
log(frequency) {\it versus} separation -- is well-described by 
the following function:
\begin{equation}
- m ( \pi_1 ) = - \frac{( 1.321 \pm .008 )}{ \log(\pi_1) } \, .
\label{mfunc}
\end{equation}
We refer the interested readers to {\it TPI} for detailed discussion of the
derivation of this result and reinforcement of the claim that this
form for $m$ is consistent with the Prime Number Theorem, the Twin
Prime Conjecture, and the Hardy--Littlewood Conjecture for twins.

\section{Two Results}

Armed with this model for the distribution of twins we now undertake
some analyses of its predictions.

\subsection{Gaps}

Those who are actively enumerating twin primes are quite interested in
the occurrences of large gaps~\cite{brent1,brent2,nicely,wolf}.
In our probabilistic model we can make estimates for the threshold
number at which we might expect particular large prime separations 
to appear.
The analysis below does just that.
There is a caveat however:  the gaps that we determine are 
prime separations, while the gaps referred to in the literature
are gaps in the natural number sequence ({\it i.e.,} arithmetic
differences between an element of a twin and the corresponding 
element of its immediate predecessor).
 
To determine likelihood thresholds for large
prime separations, we need to convolve the behaviour of $m$
as a function of $\log( \pi_1 )$ with its role in expressing
the probability (or frequency) of separations ${\cal P}$.
Simply put, the probability that a particular large prime 
separation $s_L$ occurs is (we set $ds = \Delta s = 1$)
\begin{equation}
{\cal P}(s_L, \pi_1) = \frac{m_0}{\log(\pi_1)} \, 
\exp\left(- \frac{m_0 \, s_L}{\log(\pi_1)}\right) \, .
\label{sLprob}
\end{equation}
In this expression for the probability, we have merely taken 
(\ref{prob}) and inserted our empirical fit for $m$ (\ref{mfunc})
where, for convenience, we have written $m_0$ instead of the 
central value of the constant 1.321.
The probability is the expected normalised frequency of each
type of gap among the $\pi_2$ twins (actually there are $\pi_2 -1$
gaps between the twins, but we drop the -1 in this analysis).

While this is completely correct as it stands, and provides the 
probability that a prime separation of magnitude $s_L$ occurs
somewhere within the subset of first primes with length $\pi_1$,
this not quite what is wanted.
Instead of looking for a gap of length $s_L$, we should be considering
{\it any} gap of length greater than or equal to $s_L$.
To accomplish this, we integrate ${\cal P}$ from $s_L$ to $\infty$,
and obtain the probability that a prime separation $\ge s_L$ occurs
\begin{equation}
{\cal P}_{\ge}(s_L, \pi_1) = 
\exp\left(- {\frac{m_0 \, s_L}{\log(\pi_{1})}}\right) \, .
\label{sL+prob}
\end{equation}
Let us set a {\it minimum threshold} of $f$ separation ``events''
out of the total number, $\pi_2$ , where $f$ is a risk factor 
allowing one to be more cautious ($f > 1$) or more daring ($f < 1$).
Fixing 
\begin{equation}
{\cal P}_{\ge}(s_L, \pi_1) = \frac{f}{\pi_2} \, ,
\label{sL+frac}
\end{equation}
leads to a relation among the size of the gap, the prime and 
twin counts, and the risk factor:
\begin{equation}
s_L = \frac{1}{m_0} \log(\pi_{1})  \, 
\big( \log( \pi_{2} ) - \log(f) \big) \, .
\label{sL+prob1}
\end{equation}
To cast this in a more useful form, we will utilise the well-known
approximate formulae for $\pi_1$ and $\pi_2$ quoted in {\it TPI}.
\begin{equation}
\widetilde{\pi}_{1}(N) \sim \frac{N}{\log(N)} \quad , \quad
\widetilde{\pi}_{2}(N) \sim 2 c_2 \, \frac{N}{(\log(N))^2} \, .
\label{draconian}
\end{equation}
With these approximations, (\ref{sL+prob1}) becomes
\begin{equation}
s_L = \frac{1}{m_0} \Big[ \log(N) - \log( \log(N) ) \Big] \, 
\Big[ \log(N) - 2 \log( \log(N) ) + \log( 2 c_2 ) - \log(f) \Big] \ .
\label{gapper}
\end{equation}
While this equation does not admit simple inversion to determine 
$N ( s_L )$, it is readily graphed.
We do this below 
up to prime separations of about $250$ allowing us to compare
the predictions of our probabilistic model with the data.

{\begin{figure}[htb]
\begin{center}
\begin{turn}{-90}
\leavevmode
\epsfxsize=3.5in
\epsfysize=4.5in
\epsfbox{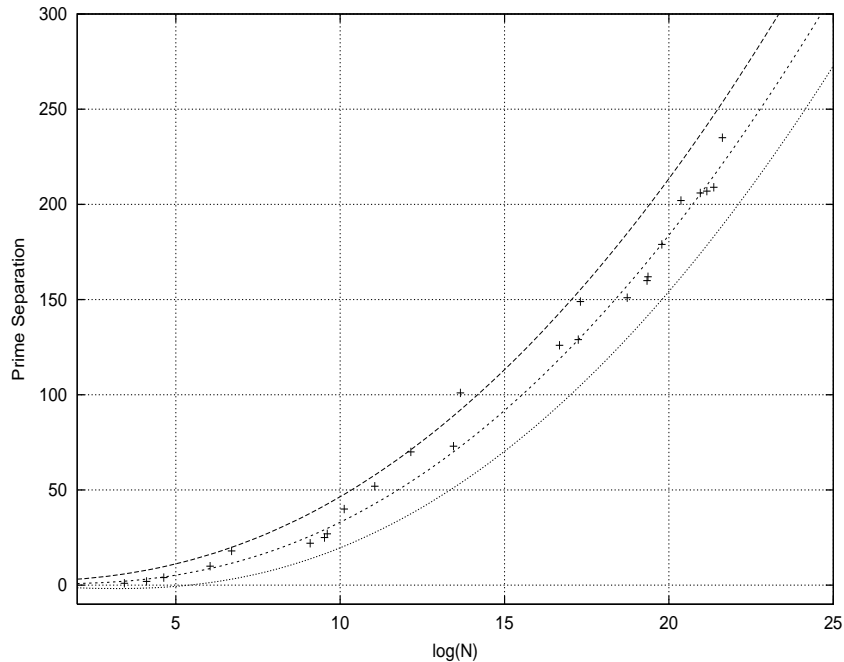}
\end{turn}
\end{center}
\caption{Prime gap {\it vs.} $\log(N)$, 
for $m_0 = 1.321$, and $0.1 \le f \le 10$.
The points mark the onset of successive maxima in the spectrum of gaps.}
\label{graph1}
\end{figure}}
 
As one can see from the graph, the probabilistic model yields 
predictions which conform extremely well with the available data.
The shape of the curves which bound our cautious--daring range
is consistent with the Hardy--Littlewood Conjecture for twins
in that the twins become more sparse among the primes as $N$ increases.

\subsection{High-Jumpers}

The high-jumper at any fixed $N$ is the natural number separation 
which occurs most frequently between consecutive primes less than $N$.
The high-jumpers are observed to increase with $N$~\cite{twinjump}.
Exactly the same analysis can be carried over to the case of twins.

We refine the notion of a twin high-jumper to 
conform to our viewpoint that it is the distribution of twins among
the primes themselves, not the natural numbers, which is significant.
In other words our twin high-jumper is the most probable prime separation
at any given $N$.
With our probabilistic decay model for the occurrence of twins, it 
necessarily follows that the {\it most likely}
high-jumper for twins is always zero.

We note that this claim is crucially dependent upon our model for the
slope $- m$ tending to $0^{-}$ (or any finite negative value) rather
than a positive value as $N \longrightarrow \infty$.
Furthermore, we recognise that since the slowly-varying probability
model for the occurrence of twins is stochastic, there may exist
ranges of $N$ for which the high-jumper is actually $1$ or $2$ or greater.

\section{Conclusion}

In this paper, we report on two results that follow
straightforwardly from our novel characterisation of the
sequence of twin primes as slowly-varying probability random events
among the primes themselves.
The first is that the model can be used to make probabilistic
predictions for the onset of appearance of prime separations of 
a given (large) size.
We have done this and the predictions agree very well 
in the regime where we have data.
We believe that the predictions extrapolate and remain valid beyond
the range of our data.
The second result is that the so-called {\it twin prime high-jumper,}
the mode prime separation, is always predicted to be zero by 
the very structure of the empirical model.

\section{Acknowledgements}

PFK and TP thank J.~Calvo and J.~Coykendall for helpful comments.
This work was supported in part by the National Science Foundation
(USA) under grant \#OSR-9452892 and a Doctoral Dissertation 
Fellowship.

\end{document}